\pdfoutput=1
\documentclass[paper=letter, fontsize=10pt,leqno]{scrartcl}
\usepackage[body={5.35in,7.5in}]{geometry} 

\usepackage[english]{babel}				
\usepackage{amsmath,tabu}				
\usepackage{amsfonts}
\usepackage{amsthm}
\usepackage{latexsym}
\usepackage{hyperref}

 \usepackage{verbatim}
\usepackage[pdftex]{graphicx}				
\usepackage{url}
\usepackage{amssymb}
\usepackage{bbm}
\usepackage{MnSymbol}
\usepackage{booktabs}
\usepackage[pdftex]{color}

\usepackage{amscd}
\usepackage{calrsfs}
\usepackage[format=plain]{caption}
\usepackage{wrapfig}
 \usepackage{float}
 
 \usepackage{bigints}

\setkomafont{sectioning}{\normalfont\scshape}

\usepackage{mathtools}
\DeclareCaptionStyle{ruled}{labelfont=bf,labelsep=space}
\usepackage[]{forest}

\newtheorem*{proposition*}{Proposition}

\title{
		\vspace{-1in} 	
		\usefont{OT1}{bch}{b}{n}
		\normalfont \normalsize \textsc{} \\ [25pt]
		\huge  Pulse-response analysis of a simple reaction-advection-diffusion equation
}

\author{\normalfont \large 
Jiasong Zhu\footnote{\scriptsize Department of Mathematics, Washington University, Campus Box 1146, St. Louis, MO 63130, USA},
\ Renato Feres\footnotemark[1],
\ Donsub Rim\footnotemark[1],
\ Gregory Yablonsky\footnote{\scriptsize Department of Energy, Environmental \& Chemical Engineering, McKelvey School of Engineering, Washington University, St. Louis, Mo 63130, USA}
}

\begin{document}

\maketitle

\begin{abstract}
\begin{center}
 Abstract \end{center}
{\small  
We analyze a reaction-advection-diffusion (RAD) equation arising in pulse-response studies of transport and reaction in a narrow reactor tube. A short pulse of gas is injected near one end of the reactor, while unreacted gas and reaction products are detected at the outlet. Particular attention is given to the effect of a constant axial advection velocity as a minimal extension of the standard diffusion model.

For a single pulse, we obtain analytical expressions for the exit-flow response and for experimentally accessible characteristics, including its moments and peak properties, as functions of the P\'eclet number and the reaction rate. The corresponding diffusion-advection model without reaction defines a standard transport curve that may be used as a baseline for identifying chemical activity. For a first-order irreversible reaction, the reactive and nonreactive exit-flow curves satisfy a simple exponential factorization, allowing the reaction rate to be extracted directly from their ratio.

We also extend the single-pulse analysis to a uniform periodic train of pulses. Using Poisson summation, we derive explicit expressions for the Fourier coefficients of the asymptotic periodic exit flow and relate them to the Laplace transform and moments of the single-pulse response. Finally, a stochastic interpretation in terms of reflected diffusion, first-passage times, and exponential killing provides probabilistic meaning for several of the analytical results and unifies the single-pulse, periodic-response, and reaction aspects of the model.

 }
\end{abstract}

\section{Introduction}

Detailed understanding of reaction mechanisms and chemical kinetics is a necessary condition for the development and optimization of catalysts and catalytic processes. While steady-state investigations are known to give a global view of a catalytic system, transient studies are invaluable because they can provide more detailed insight into the elementary steps involved in the process. Pulse-response studies have found wide applications in chemical kinetics and chemical engineering. One of the most widely used methods in such studies is the Temporal Analysis of Products (TAP) \cite{CYMGI,CYMGII}. For almost forty years, TAP has been successfully employed in transient studies of gas-phase heterogeneous reactions, with substantial advances in instrumentation and numerical modeling during this period. The TAP technique and the laboratory hardware used in it were pioneered by John Gleaves, and the general methodology for TAP data analysis was developed in a number of works, among which we cite \cite{CYMGI,CYMGII,CYMGIII,CSYMG,GYPS,MMF,MYC,YBMC,YCM,YKPG,YSCG}. See also Chapter 10 of \cite{YBGE}.

\begin{figure}[htbp]
\begin{center}
\includegraphics[width=5.0in]{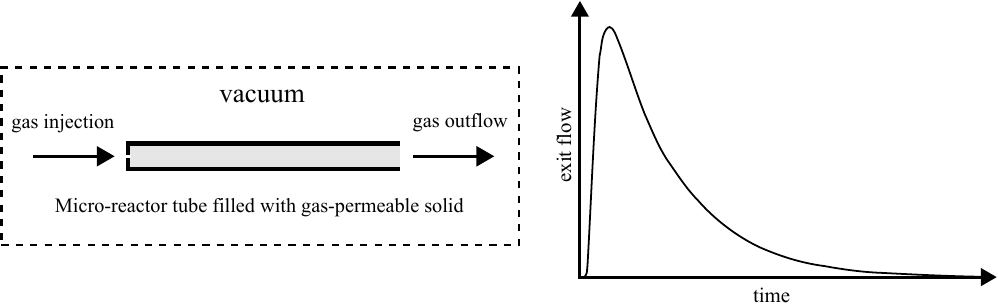}\ \
\caption{{\small Rough scheme of a TAP experiment. Small and short pulses of gas are injected into the micro-reactor tube, which is packed with solid material that is both permeable to gas diffusion and contains catalytic material promoting gas reaction. The exit gas mixture (amount of escaped gas per unit time for each gas species) is analyzed by a mass spectrometer, and the time profiles of the individual gas species (on the right) are recorded.}}
\label{TAP}
\end{center}
\end{figure}

At the core of the TAP apparatus is a micro-reactor tube of relatively small diameter filled with a catalytic sample. The sample may be embedded in a porous solid material that allows gases to diffuse through the interstitial voids and react. Small and narrow gas pulses are injected into the micro-reactor at one end of the tube, and a mixture of unconverted and newly formed gases escapes at the opposite end---the reactor outlet---where it can be analyzed by a mass spectrometer. Under typical TAP operating conditions, gas transport is well described by Knudsen diffusion, and TAP pulse-response experiments can be effectively modeled by a one-dimensional reaction-diffusion equation; see, for example, \cite{YKPG}. Figure \ref{TAP} gives a schematic description of a TAP experiment.

Under ideal TAP Knudsen-diffusion conditions, the pulse intensity is sufficiently small that intermolecular collisions are negligible. In this regime the effective diffusivity is independent of gas concentration, and the normalized response of an inert gas is independent of the amount injected. As the pulse intensity increases, however, the local pressure produced by the pulse and the frequency of gas--gas collisions also increase. Transport then enters a crossover regime in which Knudsen diffusion coexists with pressure-driven, or viscid, motion. With still larger pulses, the simple linear diffusion model is expected to break down and more complicated, nonlinear transport mechanisms become important.

More complete descriptions of transport near and beyond the Knudsen regime, such as dusty-gas and detailed pulsed-reactor models, generally involve concentration- and pressure-dependent transport coefficients and velocity fields that vary in space and time \cite{MME67}. Of particular relevance to the present work, Gering, Baroi, and Fushimi \cite{GBF18} developed a detailed model of pulsed-reactor transport covering both Knudsen and beyond-Knudsen conditions. Their model explicitly incorporates diffusive and advective contributions and describes the increasing importance of advective transport as the pulse magnitude is increased. More recently, Kunz, Maiti, Yablonsky, and Fushimi \cite{KMYF25} introduced, through a statistical analysis of experimental pulse-response data, a quantity they call the {\em rate-of-transport parameter} $\lambda$. In the Knudsen regime this parameter was found to remain essentially constant, whereas beyond that regime it exhibited a systematic dependence on a measure of the average gas concentration. These results provide further motivation for seeking simple extensions of the standard Knudsen-diffusion model that can describe the onset of concentration-dependent transport.

Motivated by these observations, we consider here a minimal linear modification of the usual transport-reaction model. Our standard transport model contains, in addition to diffusion, a constant axial advection term. The constant velocity is not intended as a complete constitutive description of transport beyond the Knudsen regime; rather, it is introduced as a simple effective parameter intended to capture, over a limited range of pulse intensities, the onset of directed pressure-driven transport. In this sense, the diffusion-advection model provides an analytically tractable extension of the standard Knudsen transport curve. The consequences and limitations of this approximation are discussed further in the body of the paper.

The introduction of an advective contribution also modifies the usual notion of a standard transport curve. A central idea in pulse-response analysis is that transport should first be characterized in the absence of reaction. Once this baseline transport response has been established from inert-gas experiments, deviations observed under reactive conditions can be used, within the assumptions of the model, to identify and characterize chemical activity. In the present work the standard transport curve therefore includes both diffusion and advection but no reaction. We illustrate the resulting transport-reaction analysis for a single first-order irreversible reaction. The same approach may serve as a template for more general reaction mechanisms.

The idea of introducing a constant advection velocity into a reaction-diffusion equation describing an inherently transient process raises several mathematical questions. We therefore study in detail the corresponding initial-boundary-value problem with a relatively small constant advection term, restricting attention primarily to the range $0\leq\mathrm{Pe}\leq 5$. Our goal is to obtain experimentally accessible characteristics of the exit flow, including its moments and peak properties, and to determine explicitly how these quantities depend on the P\'eclet number and the reaction constant. The diffusion-advection-reaction model studied here should thus be regarded as a minimal linear, moment-testable approximation to the onset of directed transport.

A second question concerns the pulse procedure itself. TAP instruments allow fine control of the timing of gas injection and can generate precisely spaced trains of small pulses in addition to isolated single pulses. It is therefore natural to ask how the response to a periodic train of pulses is related to the response to a single pulse. In particular, as the interval between pulses is varied, one passes from a regime in which successive responses are essentially isolated to one in which they overlap strongly and produce a nearly steady periodic output. This also gives a precise framework in which to examine the relation between a train of small pulses and the response generated by a larger total amount of injected gas.

We develop this connection by considering a uniform train of short pulses separated by a fixed time interval $T$. After transients have decayed, the exit flow approaches a periodic response. The natural experimentally accessible characteristics of this periodic output are its Fourier coefficients. We show that these coefficients can be obtained explicitly from the single-pulse response through the Poisson summation formula. In this way, the single-pulse and periodic-pulse experiments are treated within one mathematical framework rather than as separate models.

Finally, our analysis of the reaction-advection-diffusion initial-boundary-value problem, which is based primarily on separation of variables and eigenfunction expansions, is supplemented by a stochastic-calculus interpretation. The underlying advection-diffusion equation is the forward equation for a reflected diffusion process, and the outlet response can be interpreted in terms of first-passage times. This viewpoint provides probabilistic interpretations of the exit-flow moments, of the separation between transport and first-order reaction, and of the transfer function that appears naturally in the periodic-pulse analysis.

This paper is organized as follows. In Section \ref{sec01} we introduce the mathematical model in the form of an initial-boundary-value problem for the reaction-advection-diffusion (RAD) equation, express it in non-dimensional form, and derive a series solution for the gas concentration, making explicit the dependence on the P\'eclet number and reaction constant. In Section \ref{sec02} we investigate series convergence for both gas concentration and exit flow and show, in a numerical example, that only two terms of the series already provide a useful approximation of the exit flow. We also examine the sensitivity of the exit flow to the detailed shape of the initial gas pulse.

In the short Section \ref{sec03} we highlight the factorization
\[
j_k(L,t)=e^{-kt}j_0(L,t),
\]
from which the reaction constant $k\geq0$ can be determined from the standard transport curve $j_0(L,t)$ and the corresponding curve in the presence of reaction. In Section \ref{sec04} we derive expressions and obtain numerical values for two types of exit-flow characteristics: the {\em peak number}, defined as the product of the maximum exit flow and the time at which this maximum occurs, and the raw moments of the exit flow. In Section \ref{sec05} we discuss characteristic time scales of the transport process.

In Section \ref{section:train_of_pulses} we pass from the single-pulse problem to the response generated by a periodic train of short, weak pulses. We derive the limiting periodic exit flow and its Fourier coefficients and relate these coefficients to the single-pulse transfer function and moments. Finally, Section \ref{stochastics} develops a stochastic-calculus perspective on the problem. Although this section is principally intended to provide detailed justification and interpretation of several results obtained earlier, it also gives a unified probabilistic picture of the transport, reaction, and periodic-response aspects of the model.

  \section{The RAD initial boundary-value problem}\label{sec01}
The pulse-response system will be modeled by the following one-dimensional reaction advection diffusion equation on the interval $[0,L]$,  where $L$ is the length of the micro-reactor,
\begin{equation}\label{RAD_equation}\frac{\partial c}{\partial t}=D\frac{\partial^2 c }{\partial x^2} - v \frac{\partial c}{\partial x} - kc. \end{equation}
Here $c(x,t)$ is gas concentration for a gas species $A$, $D$ is the constant diffusion coefficient for $A$ (which depends on the characteristics of the permeable solid in the reactor's packed bed), $k\geq 0$ is the reaction constant 
for the first-order irreversible reaction $A\rightarrow B$, where $B$ is a second gas species, and $v\geq 0$ is constant advection speed. The quantity
\begin{equation}j(x,t) =c(x,t)v-D\frac{\partial c}{\partial x}(x,t) \end{equation}
is the {\em flux} of $A$.    The above model is justified in TAP experiments under a number of assumptions, among which we highlight:
   (1) the void fraction of the cylindrical  reactor's packed bed is uniform; (2) radial concentration gradient can be neglected; (3)
 axial or radial temperature  gradients can be neglected;  (4) surface diffusion can be neglected. See \cite{YKPG}.

What can be experimentally obtained in a pulse-response study is the outflow $j(L,t)$, shown schematically on the right of Figure \ref{TAP}. The  Dirichlet condition $c(L,t)=0$ will be assumed, corresponding to the  vacuum   at the right-end exit, and  no-flux/Robin condition,  $j(0,t)=0$, at the left end.
 The latter condition means that no amount of gas can enter or leave the reactor at $x=0$ for $t>0$ (after the initial injection).
  To model the initial pulse of gas, we suppose that 
  \begin{equation}c(x,0)=a\delta_\ell(x|x_0),\end{equation}
where  $\delta_\ell(x|x_0)$ is a differentiable unit mass density supported on an interval of length $\ell$ and centered around $x_0>0$ while $a$ represents the amount of gas injected. We will later pass to the limit $\ell\rightarrow 0$ and
  take, for concreteness,  $a=1$. Note that $\delta_\ell(x|x_0)$ approximates a Dirac delta function at $x_0$.

  \subsection{Non-dimensionalization of the RAD equation}
 It is convenient to rewrite the RAD equation using non-dimensional variables. Following \cite{LL}, we introduce the non-dimensional position $\xi=x/L$
 and, for $v>0$, the time variable $\tau=vt/L$, in terms of which the RAD equation becomes
 \begin{equation}\frac{\partial \rho}{\partial \tau}= \frac1{\text{Pe}} \frac{\partial^2\rho}{\partial \xi^2} -\frac{\partial \rho}{\partial \xi} -\kappa \rho, \end{equation}
 where $\rho(\xi,\tau)=c(x(\xi), t(\tau))$ is defined for $\xi$ in $[0,1]$ and 
 \begin{equation}\text{Pe}=\frac{vL}{D}, \ \ \kappa=\frac{k L}{v}.\end{equation}
 These constants are known, respectively, as the {\em P\'eclet} and {\em Damk\"ohler} numbers. The variables $\tau$ and $\kappa$ are used here for $v>0$; the case $v=0$ is treated directly, or equivalently obtained from the final dimensional formulas by taking the limit $\text{Pe}\to 0$.    The flux $j$ becomes
 \begin{equation}j(x(\xi),t(\tau))= v\left(\rho(\xi,\tau)-\frac1{\text{Pe}}\frac{\partial \rho}{\partial\xi}\right). \end{equation}
 We define the rescaled flux  $\widetilde J(\xi,\tau)=j(x(\xi),t(\tau))/v$.
 For non-zero $v$, the zero-flux condition on the left-end is the Robin condition
 \begin{equation}\rho(0,\tau)-\frac{1}{\text{Pe}}\frac{\partial \rho}{\partial \xi}(0,\tau)=0. \end{equation}
 If $v=0$, this boundary condition reduces to standard Neumann condition $\frac{\partial \rho}{\partial \xi}(0,\tau)=0$. 
 On the right-end, $\rho(1, \tau)=0$. 
Let us further introduce the function $Y(\xi,\tau)$ by the expression
\begin{equation} \rho(\xi,\tau)=\exp\left\{\frac{\text{Pe}}{2}\xi -\left(\frac{\text{Pe}}{4}+\kappa\right)\tau\right\}Y(\xi,\tau).\end{equation}
A straightforward computation yields the transformed equation
\begin{equation}\label{heat} \frac{\partial Y}{\partial \tau}= \frac1{\text{Pe}}\frac{\partial^2 Y}{\partial\xi^2}\end{equation}
and the rescaled flux, 
\begin{equation} \widetilde J(\xi,\tau)=\exp\left\{\frac{\text{Pe}}{2}\xi -\left(\frac{\text{Pe}}{4}+\kappa\right)\tau\right\} \left(\frac12 Y(\xi,\tau)-\frac1{\text{Pe}}\frac{\partial Y}{\partial \xi}(\xi,\tau)\right).\end{equation}
Boundary conditions become
\begin{equation}\label{boundary_conditions_Y}Y(1,\tau)=0, \ \ \frac12 Y(0,\tau)-\frac1{\text{Pe}}\frac{\partial Y}{\partial \xi}(0,\tau)=0\end{equation}
and the initial condition expressed for  $Y$  is
\begin{equation}\label{initial}Y(\xi,0)= \frac{a}{L} \exp\left\{-\frac{\text{Pe}}{2}\xi\right\}\delta_\epsilon(\xi|\xi_0),\end{equation}
where $\epsilon=\ell/L$, $\xi_0=x_0/L$, and $\delta_\epsilon(\xi|\xi_0):=L\,\delta_\ell(L\xi|x_0)$ is a unit-mass density on the $\xi$-interval.
\subsection{Series solution}
As noted above, the RAD equation is reduced to the simple diffusion equation (\ref{heat}). Its solution can be obtained by  the standard separation of variables method (see, for example, \cite{Logan}).
If  we set \begin{equation}Y(\xi,\tau)=T(\tau)\Phi(\xi),\end{equation} then Equation (\ref{heat}) implies 
\begin{equation}{\text{Pe}}  \frac{T'(\tau)}{T(\tau)} =\frac{\Phi''(\xi)}{\Phi(\xi)}=-\mu^2.\end{equation}
Note the choice of negative constant. This is because solutions should decay to zero for large $\tau$ as the reactor eventually fully evacuates. Thus, up to a multiplicative constant, 
\begin{equation} T(\tau)=\exp\left(-\frac{\mu^2}{\text{Pe}}\tau\right)\end{equation}
and
\begin{equation}\label{Phi_xi}
\Phi(\xi)=A\mu\cos(\mu\xi)+B\sin(\mu\xi).
\end{equation}
The boundary conditions of Equation (\ref{boundary_conditions_Y}) imply for  $\Phi(\xi)$:
\begin{equation}\label{boundary_conditions_Phi}
\Phi(1)=0, \ \ \frac12\Phi(0)-\frac1{\text{Pe}}\Phi'(0)=0, 
\end{equation}
from which we obtain that $\mu>0$ must satisfy the equation
\begin{equation}\label{eigenvalues}
\cos\mu =-\frac{\text{Pe}}2 \frac{\sin\mu}{\mu}
\end{equation}
and, up to a multiplicative constant,  
\begin{equation}\label{Phi_xi_bis}
\Phi(\xi)={\sin\left(\mu(1-\xi)\right)}.
\end{equation}
The solutions to Equation (\ref{eigenvalues})  form an increasing sequence $\mu_n$ which, for large $n$, has the asymptotic expression
\begin{equation}
\mu_n\approx \mu^{\text{\tiny asymp}}_{n}= \left(n-\frac12\right)\pi,
\end{equation}
where $n$ is a positive integer. 
In order to obtain approximate values for $\mu_n$, we rewrite Equation (\ref{eigenvalues}) as
\begin{equation}
f(\mu)=\frac{2}{\text{Pe}}\mu + \tan \mu=0.
\end{equation}
For $\text{Pe}>0$, $f'(\mu)>0$ on each interval
\begin{equation}I_n=\left(\left(n-\frac12\right)\pi,n\pi\right), \qquad n=1,2,\ldots .\end{equation}
Moreover, $f(\mu)\to-\infty$ as $\mu$ approaches the left endpoint of $I_n$ from the right, while $f(n\pi)>0$. Therefore $f$ contains a single zero, $\mu_n$, in each $I_n$, which can
 be efficiently  approximated by the standard bisection method with good control on error bounds. For example,
 to ensure that the bracketing interval has length at most $10^{-8}$, it is enough to take
 $m\geq \log_2\left(10^8\pi/2\right)$, hence $m\geq 28$. The case $\text{Pe}=0$ is obtained directly from Equation (\ref{eigenvalues}), giving $\mu_n=(n-\frac12)\pi$. 
  See Table (\ref{tab:examplePe4}) for the values of $\mu_n$ up to $n=14$, in which  
  we set   $\text{Pe}=4$.   

\begin{table}[h!]
    \centering
    \caption{Eigenvalues $\mu_n$ ($n\leq 14$) for $\text{Pe}=4$}
    \label{tab:examplePe4}
    \begin{tabular}{c| c c c c c c c c}
      \toprule n & 1 & 2 & 3 & 4 & 5 & 6 & 7 \\ $\mu_n$ & 2.2889 & 5.0870 & 8.0962 & 11.1727 & 14.2764 & 17.3932 & 20.5175 \\\hline n & 8 & 9 & 10 & 11 & 12 & 13 & 14 \\ $\mu_n$ & 23.6463 & 26.7781 & 29.9119 & 33.0472 & 36.1835 & 39.3207 & 42.4586\\\bottomrule   \end{tabular}
\end{table}

\begin{figure}[htbp]
\begin{center}
\includegraphics[width=3.5in]{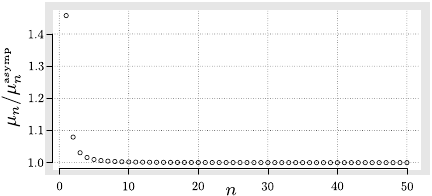}\ \ 
\caption{{\small  Ratio of  $\mu_n$ over the asymptotic value $\mu^{\text{\tiny asymp}}_n=\left(n-\frac12\right)\pi$.
}}
\label{compare}
\end{center}
\end{figure}

For the same value of the P\'eclet parameter ($\text{Pe}=4$), Figure \ref{compare} shows   that  $\mu^{\text{\tiny asymp}}_n$ becomes a fairly good approximation of $\mu_n$ for $n\geq 10$.

The functions
\begin{equation}\label{Phi_n_xi}\Phi_n(\xi)= \frac{\sin(\mu_n(1-\xi))}{\sqrt{\frac12+\frac{\text{Pe}}{4}\left(\frac{\sin\mu_n}{\mu_n}\right)^2}}\end{equation}
are the normalized eigenfunctions of the self-adjoint Sturm--Liouville operator $-d^2/d\xi^2$ with the boundary conditions (\ref{boundary_conditions_Phi}), and thus comprise an orthonormal basis for $L^2([0,1])$.  See
\cite{Strauss}.

We can then search for a solution of the initial boundary-value problem by setting
\begin{equation}
Y(\xi,\tau)=\sum_{n=1}^\infty A_n\exp\left(-\frac{\mu_n^2}{\text{Pe}}\tau\right) \Phi_n(\xi)
\end{equation}
or, back to $\rho(\xi,\tau)$,
\begin{equation}
\rho(\xi,\tau)=\exp\left\{\frac{\text{Pe}}{2}\xi-\left(\frac{\text{Pe}}{4}+\kappa\right)\tau\right\} \sum_{n=1}^\infty A_n\exp\left(-\frac{\mu_n^2}{\text{Pe}}\tau\right) \Phi_n(\xi).
\end{equation}
The coefficients $A_n$ are obtained by imposing the initial condition (\ref{initial}), or
\begin{equation}
\rho(\xi,0)=\frac{a}{L}\delta_\epsilon(\xi|\xi_0).
\end{equation} 
Therefore
\begin{equation}
A_n=\frac{a}{L}\int_0^1 \delta_\epsilon(\xi|\xi_0)e^{-\frac{\text{Pe}\xi}{2}}\Phi_n(\xi)\, d\xi.
\end{equation}

The ideal short   pulse of gas corresponds to the limit for $\epsilon$ approaching $0$.  Taking this limit,   we obtain
\begin{equation}
A_n=\frac{a}{L} e^{-\frac{\text{Pe}\xi_0}{2}}\Phi_n(\xi_0)= \frac{(a/L)e^{-\frac{\text{Pe}\xi_0}{2}}\sin(\mu_n(1-\xi_0))}{\sqrt{\frac12+\frac{\text{Pe}}{4}\left(\frac{\sin\mu_n}{\mu_n}\right)^2}}.
\end{equation}

The solution to the initial boundary-value problem, therefore, takes the form
\begin{equation}\label{solution_xi_tau}
\rho(\xi,\tau)=\frac{a}{L}\exp\left\{\frac{\text{Pe}}{2}(\xi-\xi_0)-\left(\frac{\text{Pe}}{4}+\kappa\right)\tau\right\} \sum_{n=1}^\infty \Phi_n(\xi_0) \exp\left(-\frac{\mu_n^2}{\text{Pe}}\tau\right)\Phi_n(\xi).
\end{equation}

\section{Convergence of the infinite series}\label{sec02} 
\begin{figure}[htbp]
\begin{center}
\includegraphics[width=4.0in]{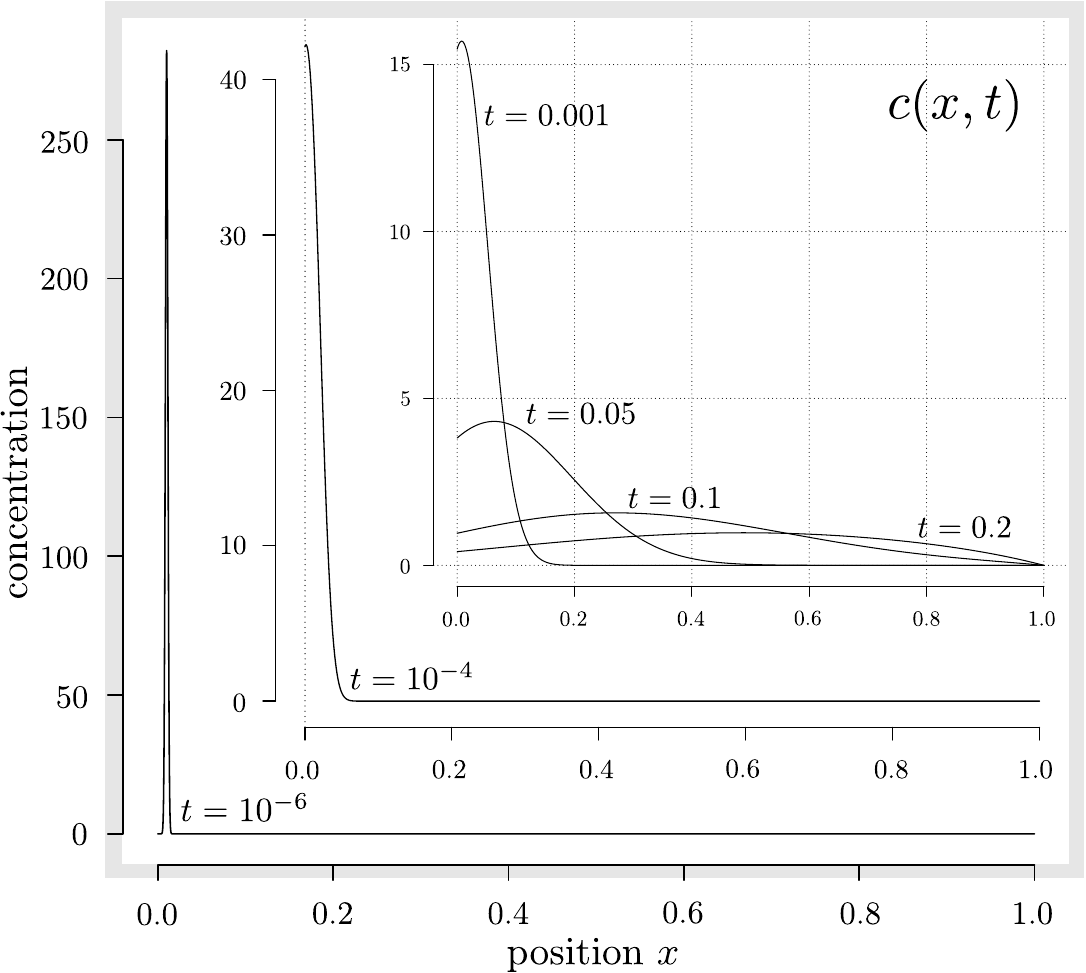}\ \ 
\caption{{\small  Profile of concentration $c(x,t)$ inside the reactor ($0\leq x\leq 1$) at different times. Parameters used: $x_0=0.01$, $k=1$, $\text{Pe}=4$. 
(Here we used $a=1$, $L=1$;  note that $k$ and $t$ are used, not $\kappa$ and $\tau$), and $n=1500$. }}
\label{concentration_time}
\end{center}
\end{figure}

For every $\tau_0>0$, the infinite series of Equation (\ref{solution_xi_tau}) converges uniformly on $[0,1]\times[\tau_0,\infty)$ and gives a smooth solution to the boundary-value problem.
Figure \ref{concentration_time} illustrates the convergence of the solution toward the delta-pulse initial condition as time decreases to zero. 

Observe that $\tau/\text{Pe}=t/t_d$, where $t_d=L^2/D$ is a characteristic time of the diffusion process (this is twice the expected time a gas molecule takes  to exit the micro-reactor when $v=0$), and $\kappa\tau=kt$. See Section \ref{sec05}.

In terms of the original variables $x$ and  $t$, the concentration function is now
\begin{equation}\label{infinite_series_c}
c(x,t)=\frac{a}{L}e^{-kt}\exp\left(\frac{\text{Pe}}{2L}(x-x_0)-\frac{\text{Pe}^2}{4t_d}t\right) \sum_{n=1}^\infty \Phi_n(x_0/L) \exp\left(-\frac{\mu_n^2}{t_d}t\right)\Phi_n(x/L).
\end{equation}
Recall that $x_0>0$ is the injection point, close to $0$ (the left-end of the micro-reactor tube).

\begin{figure}[htbp]
\begin{center}
\includegraphics[width=4.5in]{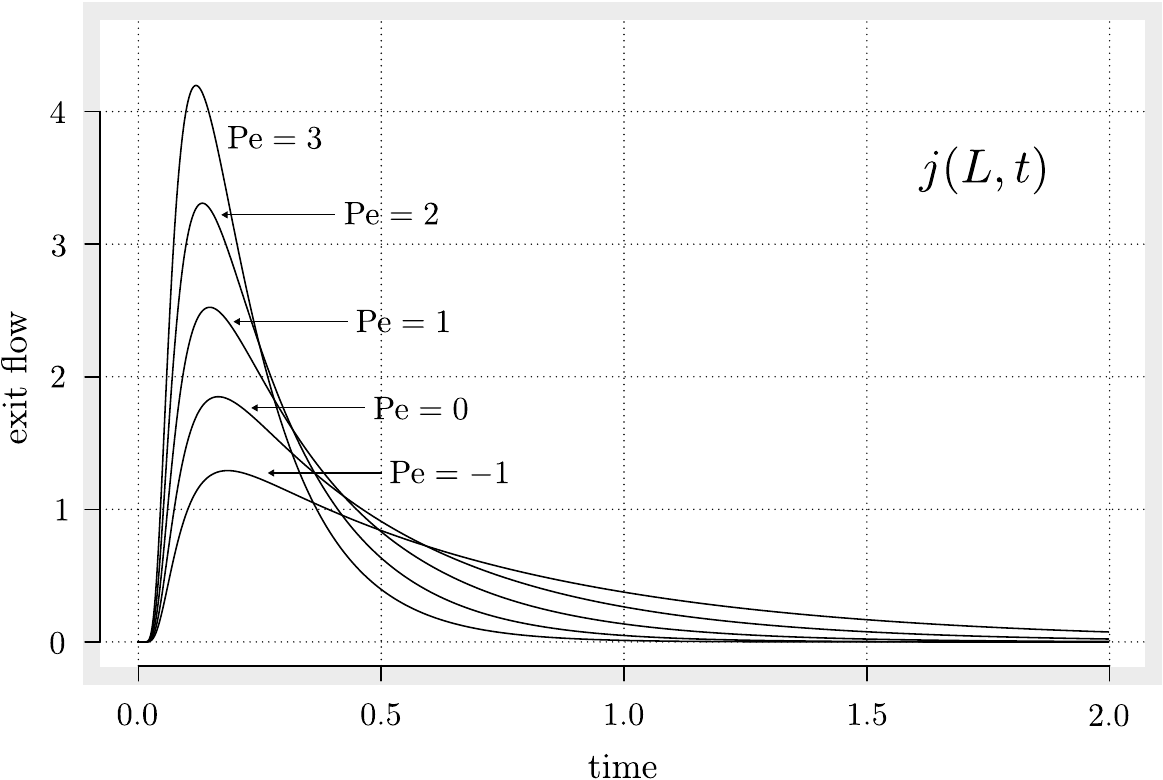}\ \ 
\caption{{\small Exit flow as a function of $\text{Pe}$. Here, $k=0$ and $x_0=0.01$. The analysis below focuses on the physically relevant range $\text{Pe}\geq0$. These are called the {\em standard transport curves}.}}
\label{Peclet}
\end{center}
\end{figure}

 The outflow of gas is the flux $j(x,t)$ at $x=L$:
 \begin{equation}
 j(L,t)=  c(L,t)v- D\frac{\partial c}{\partial x}(L,t) =  -D\frac{\partial c}{\partial x}(L,t).
 \end{equation}
Since $\Phi_n(1)=0$ and  
\begin{equation}
\Phi_n'(1)=-\frac{\mu_n}{\sqrt{\frac12+\frac{\text{Pe}}{4}\left(\frac{\sin\mu_n}{\mu_n}\right)^2}}
\end{equation}
by Equation (\ref{Phi_n_xi}), we have
 \begin{equation}\label{jL_2}
 j(L,t)=     \frac{aD}{L^2}\exp\left\{\frac{\text{Pe}}{2L}(L-x_0)-\left(\frac{\text{Pe}^2}{4t_d}+k\right)t\right\}  \sum_{n=1}^\infty \frac{\mu_n\sin\left(\mu_n\left(1-\frac{x_0}{L}\right)\right)}{\frac12+\frac{\text{Pe}}{4}\left(\frac{\sin\mu_n}{\mu_n}\right)^2} \exp\left(-\frac{\mu_n^2}{t_d}t\right).
 \end{equation}
 When convenient, we  write (\ref{jL_2}) as  $j_k(L,t)$ to make explicit the dependence on $k$. The {\em standard transport curves} are the graphs of $j_0(L,t)$.

 \subsection{Two-terms approximation of exit flow}
\begin{figure}[htbp]
\begin{center}
\includegraphics[width=2.7in]{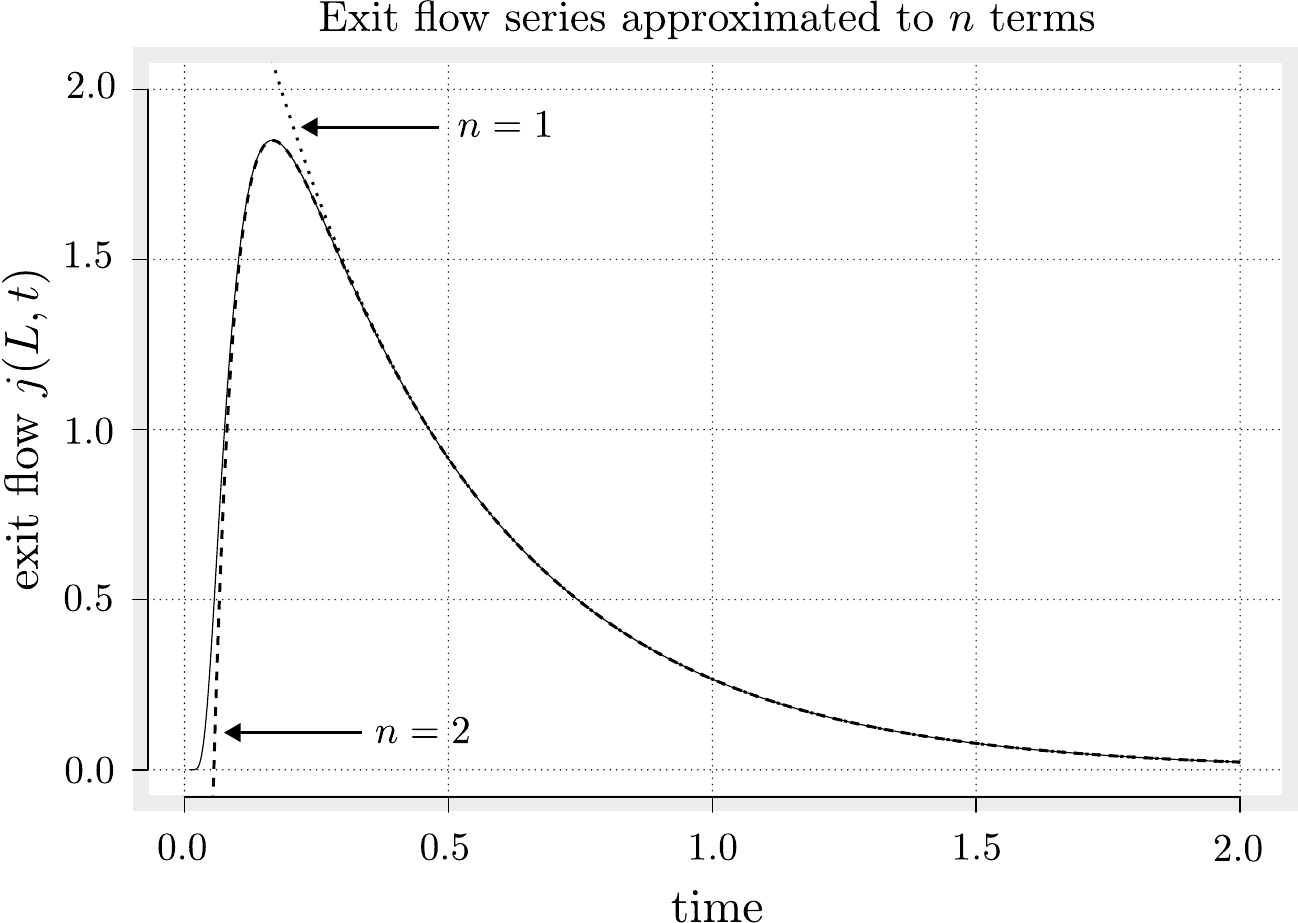}\ \ 
\caption{{\small  Approximation of $j(L,t)$ by one and two terms of the infinite series. For not too small values of $t$, two terms already provide a very useful description of the exit flow. Here $k=0$ and $\text{Pe}=0$.}}
\label{Approximation}
\end{center}
\end{figure}

 We note that the infinite series giving the exit flow in Equation (\ref{jL_2}) converges fast.  Figure \ref{Approximation}
 shows the approximation of $j(L,t)$ by the first term of the series (dotted curve), and by the first two terms (dashed curve). The thin solid curve
 uses $100$ terms.
 
\begin{figure}[htbp]
\begin{center}
\includegraphics[width=2.5in]{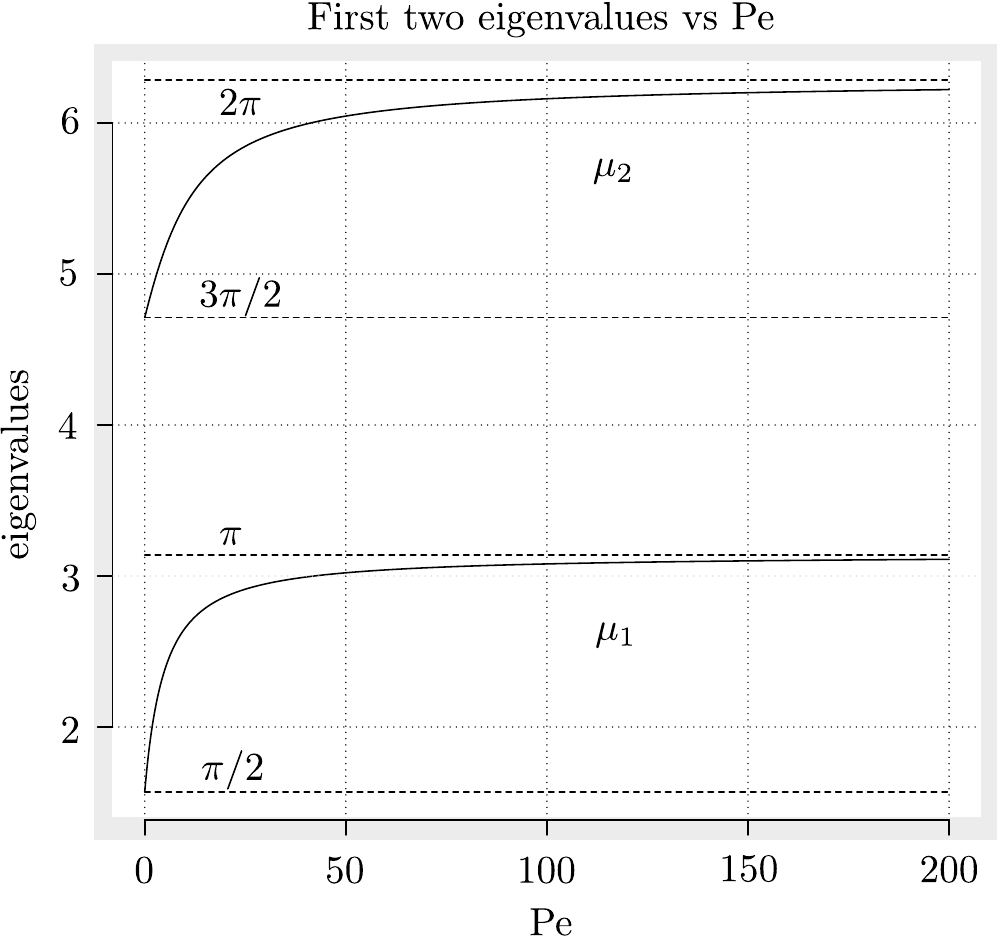}\ \ 
\caption{{\small The first two eigenvalues, $\mu_1, \mu_2$, as functions of the P\'eclet parameter. }}
\label{First_two_eig}
\end{center}
\end{figure}

 For most practical purposes, a two-term approximation of the exit flow already gives a very useful description. See, for example, the below Section \ref{section_peak}.
 The first two eigenvalues are shown in Figure \ref{First_two_eig}.
 
\subsection{Sensitivity to the shape of the  initial pulse}
 It is natural to ask how sensitive is   the exit flow to the precise form of the initial pulse. Above we have assumed that the pulse  is a delta function of the position variable $x$, supported at $x=0$. Let us compare the exit flow obtained under that assumption with the exit flow resulting from 
 a crude approximation of the delta-pulse, in the form
 \begin{equation}\delta_\epsilon(\xi|\xi_0)=\frac1\epsilon \mathbbm{1}_{[0,\epsilon]}(\xi),
 \end{equation}
 where   $\xi_0=\epsilon/2$. Here $\mathbbm{1}_{[0,\epsilon]}(\xi)$  represents a square function which is $1$ for $\xi$ in the interval $[0,\epsilon]$
 and $0$ elsewhere.
  Note that this function does not satisfy the zero flux condition at $0$. Nevertheless, it produces exit flow functions nearly indistinguishable from
  the delta-function  limit as the below table shows.
 For concreteness, we suppose $\text{Pe}=k=0$.   Letting  $j_\epsilon(t)$ denote the exit flow for this $\epsilon$-dependent initial pulse and $j_\delta(t)$ the exit flow for the
 delta-pulse, we write
 $E(\epsilon)=\|j_\epsilon-j_\delta\|$ for the maximum   of $|j_\epsilon(t)-j_\delta(t)|$ (for  $t$ in the   interval $[0,2]$,  during which  time period more than $99\%$ of the injected gas has escaped). 
 Truncating  the series expansion  of the exit flow to $100$ terms, we find the following numerically obtained values (Table \ref{tab:errors}):
 
 \begin{table}[h!]
    \centering
    \caption{ $E(\epsilon)=\max_{0\leq  t\leq  2} |j_\epsilon(t)-j_\delta(t)|$}
    \label{tab:errors}
    \begin{tabular}{c| c c c c c }
      \toprule $\epsilon$ & $10^{-1}$ & $10^{-2}$ & $10^{-3}$ & $10^{-4}$ & $10^{-5}$  \\ $E(\epsilon)$ & $0.044$ & $2\times 10^{-2}$ & $3\times 10^{-4}$& $3\times 10^{-6}$ & $3\times 10^{-8}$ \\\bottomrule   \end{tabular}
\end{table}
 The length of the micro-reactor is $1$. For the very crude approximation in which $\epsilon$ is one-tenth of that length, $E(\epsilon)$ 
 $\approx 0.05$. This maximum value occurs at relatively early times when the exit flow curve is steepest. The graphs of $j_{0.1}(t)$ and $j_\delta(t)$ are, visually, nearly indistinguishable. This is evidence that the results about exit flow obtained in this paper    are fairly insensitive to the exact shape of the initial gas pulse.

\section{Separation of reaction term}\label{sec03}

 Since the ultimate goal of pulse-response  studies is to   probe chemical kinetics properties from a baseline of pure transport,
 it is worth noting that the ratio $j_k(L,t)/j_0(L,t)$ (in which we indicate the reaction constant as a subscript) assumes the simple form
 \begin{equation}\label{ratio}
 \frac{j_k(L,t)}{j_0(L,t)}= e^{-kt}.
 \end{equation}
 This is an immediate consequence of Equation (\ref{jL_2}).
 Therefore, the reaction constant is readily obtained from a pulse-response experiment.

\section{Numerical signatures of the exit flow}\label{sec04}
Certain numerical characteristics of the exit flow that can in principle be obtained empirically from pulse-response  experimental data are worth highlighting. Here we emphasize the raw moments of the exit flow (see \cite{YSCG}) and a number associated to its peak value.
Note: our analysis is limited to $\text{Pe}\geq0.$ However, in the moment plots below we have included negative values as well. Although 
the negative range for $\text{Pe}$ is unphysical under the  assumptions of the present paper, it is mathematically meaningful,   and
possibly of interest in settings where a negative advection velocity can be imposed through mechanical means.

\subsection{$J(\tau_d) \times \tau_d$ at peak}\label{section_peak}
In order to define the exit flow signature at peak, we set $x_0=0$, normalize $j(L,t)$ by dividing it by $aD/L^2$, and make it a function of the
 non-dimensional time $\tau_d = t/t_d$:
\begin{equation}
  J(\tau_d)=     \exp\left\{\frac{\text{Pe}}{2}-\left(\frac{\text{Pe}^2}{4}+\kappa_d\right)\tau_d\right\}  \sum_{n=1}^\infty \frac{\mu_n\sin\mu_n}{\frac12+\frac{\text{Pe}}{4}\left(\frac{\sin\mu_n}{\mu_n}\right)^2} \exp\left(-{\mu_n^2}\tau_d\right).
\end{equation}
Here we have used $\kappa_d=kt_d$. 

\begin{figure}[htbp]
\begin{center}
\includegraphics[width=4.5in]{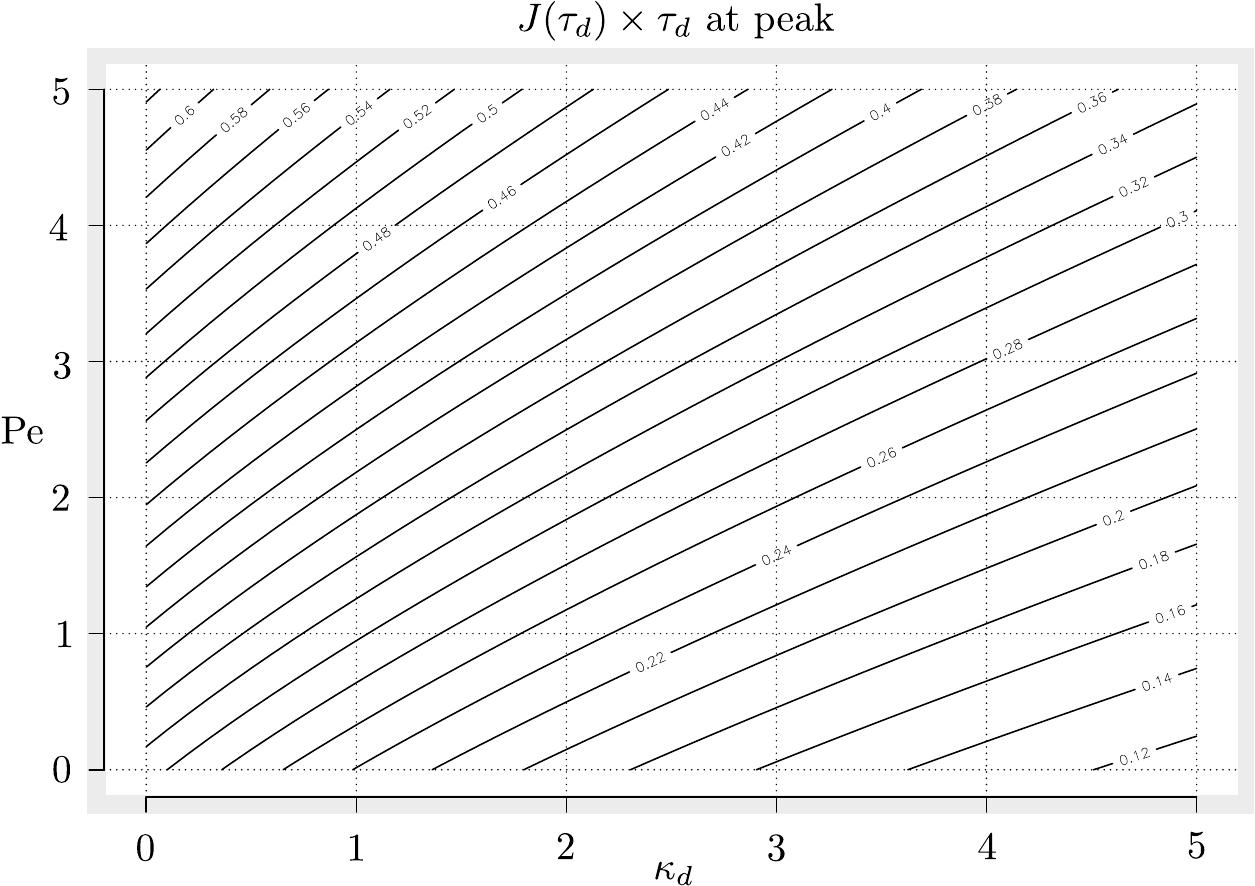}\ \ 
\caption{{\small  Product   $J(\tau_d)\times \tau_d$ at the time $\tau_d$ when the exit flow is at its maximum value as a function of the parameters $\text{Pe}$ and $\kappa_d$. Here $J(\tau_d)$ is the normalized (non-dimensional) exit flow, $\tau_d$ is the non-dimensional time $t/t_d$, and $\kappa_d=k t_d$ is non-dimensional reaction constant.  When $\text{Pe}=0, \kappa_d=0$, this signature value is $0.3083$}}
\label{Max_char}
\end{center}
\end{figure} 

 We now look for the quantity: $\tau_{\text{max}}\times \max_{\tau_d} J(\tau_d)$, where $\tau_{\text{max}}$ is the dimensionless time at which exit flow is at a maximum.
 For these parameters we have
 $$ J(\tau_d)=\pi\sum_{n=1}^\infty (-1)^{n+1}(2n-1) e^{-\left(n-\frac12\right)^2\pi^2 \tau_d}.$$
The time of maximum is the zero  of the derivative
\begin{equation}\label{Jprime} 
\begin{aligned}
J'(\tau_d)&=  -   \exp\left\{\frac{\text{Pe}}{2}-\left(\frac{\text{Pe}^2}{4}+\kappa_d\right)\tau_d\right\}\times\\
&\ \ \ \ \ \ \ \ \ \ \ \ \ \ \ \ \   \sum_{n=1}^\infty \frac{\mu_n\sin\mu_n}{\frac12+\frac{\text{Pe}}{4}\left(\frac{\sin\mu_n}{\mu_n}\right)^2} \left(\frac{\text{Pe}^2}{4} +\kappa_d +\mu_n^2\right) \exp\left(-{\mu_n^2}\tau_d\right).
\end{aligned}
\end{equation}

When $\kappa_d=\text{Pe}=0$, the product $\tau_\text{\tiny max}\times J(\tau_\text{\tiny max})$ is $0.3083$. (Cf. \cite{GYPS}.) Any experimental deviation from this value may indicate chemical activity, advective transport, or departure from the range of validity of the Knudsen-diffusion model.
The remark made in Figure \ref{Approximation} justifies using two-term approximation of $J'(\tau_d)$ to obtain the time of maximum of $J(\tau_d)$ and maximum value.
The time of maximum is
$$ \tau_\text{\tiny max} = \frac1{\mu_2^2-\mu_1^2}\log \left\{-\frac{\frac{\mu_2\sin\mu_2}{\frac12+\frac{\text{Pe}}{4}\left(\frac{\sin\mu_2}{\mu_2}\right)^2} \left(\frac{\text{Pe}^2}{4} +\kappa_d +\mu_2^2\right)}{\frac{\mu_1\sin\mu_1}{\frac12+\frac{\text{Pe}}{4}\left(\frac{\sin\mu_1}{\mu_1}\right)^2} \left(\frac{\text{Pe}^2}{4} +\kappa_d +\mu_1^2\right)}\right\},$$
where $\log$ is  natural  logarithm. For $\text{Pe}=0$, $\kappa_d=0$, we have $$\tau_\text{\tiny max} =\frac{3\log 3}{2\pi^2}\approx 0.167.$$
The exit flow at this moment (up to the first two series terms for $\text{Pe}=0$, $\kappa_d=0$) is
$$J(\tau_\text{\tiny max})=2\left(\mu_1e^{-\mu_1^2\tau_\text{\tiny max}}-\mu_2e^{-\mu_2^2\tau_\text{\tiny max} }\right) =\pi\left(e^{-\frac{\pi^2}{4}\tau_\text{\tiny max} }-3e^{-\frac{9\pi^2}{4}\tau_\text{\tiny max} }\right)\approx 1.850.$$
The product of these two values is approximately $0.309$, to be compared with the more precise $0.3083$ obtained above using a much larger number of terms.

Using the same two-term approximation, for $0\leq\text{Pe}\leq 5$ and $\kappa_d=0$, the peak characteristic $\text{Peak}(\text{Pe})=J(\tau_\text{max})\times \tau_\text{max}$ satisfies
\begin{equation}
\text{Peak}(\text{Pe})=0.308 + 0.066\times \text{Pe}\pm 0.01, \ \ \ 0\leq\text{Pe}\leq 5.
\end{equation}

\subsection{Raw moments of exit flow}
Before computing the moments of the exit flow, we note the following
  alternative way of expressing it. By an immediate application of the fundamental theorem of calculus and the assumption of zero flux at $0$ we obtain  the conservation law
\begin{equation}\label{jL_1}
j(L,t) = -\frac{d}{dt}\int_0^L c(x,t)\, dx - k \int_0^L c(x,t)\, dx.
\end{equation}

The total amount of injected gas still remaining in the reactor by time $t$ is easily obtained by term-by-term integration.  
\begin{equation}\label{Int_c}
\begin{aligned}
\int_0^L c(x,t)\, dx &= a\exp\left\{\frac{\text{Pe}}{2L}(L-x_0)-\left(k+\frac{\text{Pe}^2}{4t_d}\right)t\right\}\times\\
& \ \ \ \ \ \ \ \ \ \ \ \ \ \ \ \ \ \ \ \sum_{n=1}^\infty \frac{\sin\left(\mu_n\left(1-\frac{x_0}{L}\right)\right)e^{-\frac{\mu_n^2}{t_d}t}}{\mu_n\left[\frac12 + \frac{\text{Pe}}{4}\left(\frac{\sin\mu_n}{\mu_n}\right)^2\right]\left[1+\left(\frac{\text{Pe}}{2\mu_n}\right)^2\right]}.
\end{aligned}
\end{equation}
Applying (\ref{jL_1}) to (\ref{Int_c}) indeed gives (\ref{jL_2}), as   an elementary  calculation readily  shows. 
We may further pass to the limit $x_0\rightarrow 0$ to place the point of gas injection at $x_0=0$.

\begin{figure}[htbp]
\begin{center}
\includegraphics[width=4.5in]{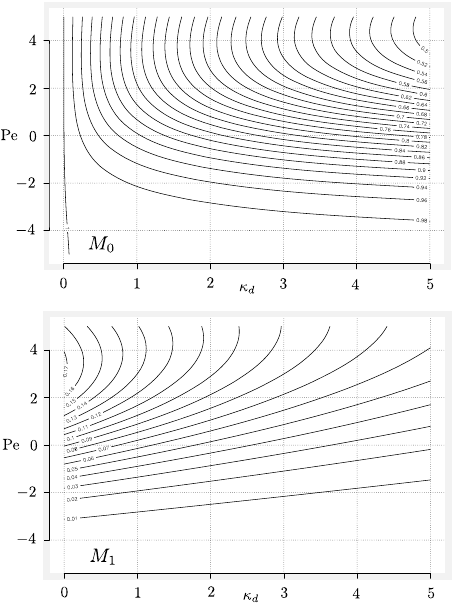}\ \ 
\caption{{\small  The first $2$ raw moments, $M_0, M_1$, of the exit flow as functions of the two parameters $\text{Pe}$ and $\kappa_d$. We have set $a=1$   and $t_d=1$. Moment $M_2$ is shown in Figure \ref{Moments2}.}}
\label{Moments}
\end{center}
\end{figure} 

The $m$th (raw)  {\em moment} of the exit flow  is defined by
$$M_m = \int_0^\infty t^m j(L,t)\, dt. $$
These are useful signature numbers of the exit flow that summarize its main properties. Our goal is to obtain values for the moments as functions of the
parameters $\text{Pe}$ and $k$. 
With this in mind, 
we define the auxiliary (non-dimensional) numbers
$$G_m=\frac{1}{t_d}\int_0^\infty\left(\frac{t}{t_d}\right)^m\int_0^L \frac{c(x,t)}{a}\, dx\, dt. $$
Using Equation (\ref{Int_c}) and the definition of the standard Gamma-function, we obtain $G_m$ as an infinite series:
\begin{equation}\label{G_m}
\begin{aligned}
G_m  &= \exp\left\{\frac{\text{Pe}}{2}\left(1-\frac{x_0}L\right)\right\}\times \\
&\ \ \ \ \ \ \ \ \ \ \  \sum_{n=1}^\infty \frac{\sin\left(\mu_n\left(1-\frac{x_0}{L}\right)\right)}{\mu_n\left[\frac12 + \frac{\text{Pe}}{4}\left(\frac{\sin\mu_n}{\mu_n}\right)^2\right]\left[1+\left(\frac{\text{Pe}}{2\mu_n}\right)^2\right]} \frac{m!}{\left(kt_d+\frac{\text{Pe}^2}{4} +{\mu_n^2}\right)^{m+1}}.
\end{aligned}
\end{equation}

A relationship between the $M_m$ and $G_m$ is obtained through Equation (\ref{jL_1}). Writing 
$I(t)=\int_0^L c(x,t)\, dx$, first observe that Equation (\ref{RAD_equation}) implies
\begin{equation}
j(L,t)  =-\frac{dI}{dt}(t) - kI(t).
\end{equation}

Note that $I(0)=a$ (where $a$ is the initial pulse intensity) and $I(\infty)=0$ (reactor fully evacuated in the long run, which happens exponentially fast).
Then
\begin{equation}
M_0 = \int_0^\infty j(L,t)\, dt =a\left(1- G_0kt_d\right).
\end{equation}
 Exponential decay of $I(t)$ implies for  $m\geq 1$ that $\int_0^\infty \frac{d}{dt}\left(t^mI(t)\right)\, dt=0$. Integration by parts and Equation (\ref{jL_1}) then yield
 \begin{equation}
 \begin{aligned}
 M_m&= - \int_0^\infty t^m I'(t)\, dt - k\int_0^\infty t^m I(t)\, dt = m\int_0^\infty t^{m-1} I(t)\, dt -k\int_0^\infty t^m I(t)\, dt\\
 &= am t_d^m  G_{m-1} - akt_d^{m+1}G_m.
 \end{aligned}
 \end{equation}

Let us record these equations here: if $m\geq 1$, and recalling that $t_d=L^2/D$, then
\begin{equation}\label{M_mG_m}
M_0=a\left(1-G_0kt_d\right), \ \ \ M_m =at_d^m \left(m G_{m-1} - k t_dG_m\right).
\end{equation}

These together with Equation (\ref{G_m}) provide  an effective way for obtaining the raw moments $M_m$.  Let us assume for concreteness that the point of gas injection is $x_0=0$. Also let us introduce the non-dimensional $\kappa_d=kt_d$ (we find this number, the second Damk\"ohler parameter, more appropriate to our needs than the velocity dependent $\kappa$). Then
$G_m=G_m(\text{Pe},\kappa_d)$, a function of only  $\kappa_d$ and the P\'eclet number.
Rewriting $G_m$ (from Equation (\ref{G_m})) to incorporate these definitions and assumptions,  
\begin{equation}\label{G_mNew}
\begin{aligned}
G_m(\text{Pe},\kappa_d) & = 2m! \exp\left(\frac{\text{Pe}}{2}\right)\times\\
 &\ \ \ \ \ 
 \sum_{n=1}^\infty 
 \frac{{\mu^{-2m-3}_n}\sin\mu_n}{\left[1+ \frac{\text{Pe}}{2}\left(\frac{\sin\mu_n}{\mu_n}\right)^2\right]\left[1+\left(\frac{\text{Pe}}{2\mu_n}\right)^2\right]\left[1+\left(\frac{\text{Pe}}{2\mu_n}\right)^2+ \frac{\kappa_d}{\mu_n^2}\right]^{m+1}}.
\end{aligned}
\end{equation}

Figure \ref{Moments} gives the first $3$ moments as functions of the parameters $\text{Pe}$ and $\kappa_d$.
\begin{figure}[htbp]
\begin{center}
\includegraphics[width=4.5in]{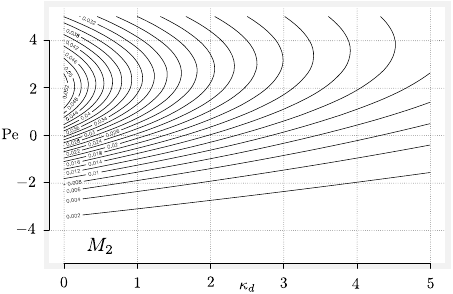}\ \ 
\caption{{\small  The third raw moments, $M_2$, of the exit flow as function of the two parameters $\text{Pe}$ and $\kappa_d$. We have set $a=1$   and $t_d=1$.}}
\label{Moments2}
\end{center}
\end{figure} 

\section{Remark about time scales}\label{sec05}
A natural time scale for the diffusion-advection transport process is set by the mean time  a gas molecule takes from the moment of injection   until it leaves the micro-reactor. Standard stochastic methods (see the last paragraph in this section) can be used to show that this mean time is given by the function $T(x)$ of the point of injection $x$
that solves the boundary-value problem
\begin{equation}\label{expected_exit}
DT''(x) +vT'(x)=-1
\end{equation}
with boundary conditions $T(L)=0$ and $T'(0)=0$. This easily  yields the expected time starting at $x=0$:
\begin{equation} 
t_{\text{\tiny mean}}= T(0)=\frac{D}{v^2}\left(e^{-\frac{vL}{D}}-1 + \frac{vL}{D}\right).
\end{equation}
In terms of the non-dimensional parameter $\text{Pe}$ and the previously defined   $t_d=L^2/D$, this is
\begin{equation}\label{mean_time}
t_{\text{\tiny mean}}=t_d \frac{ e^{-\text{Pe}}-1 + \text{Pe}}{\text{Pe}^2}.
\end{equation}

When $\text{Pe}$ is small (non-negative), diffusion dominates, and we obtain the diffusion time scale
\begin{equation} t_{\text{\tiny mean}}\approx\frac{t_d}{2}=\frac{L^2}{2D}.\end{equation}
When $\text{Pe}$ is large, advection dominates, and we obtain the advection time scale
\begin{equation}
t_{\text{\tiny mean}}\approx \frac{L}{v},
\end{equation}
which is known as the {\em residence time} in flow systems. Yet another time-related quantity of direct empirical access is the ratio of moments
\begin{equation}
t_{\text{\tiny moments}}=\frac{M_1}{M_0}.
\end{equation}
For the inert case $k=0$, it turns out that 
\[t_{\text{\tiny moments}}=t_{\text{\tiny mean}}.\] See section \ref{stochastics}.

At the beginning of this section we referred to stochastic methods. Briefly, the advection-diffusion equation  is the Fokker-Planck equation associated to
the stochastic differential equation $dX_t = vdt +\sqrt{2D} dW_t$ that models one-dimensional molecular motion. An application of Dynkin's formula from the theory of stochastic differential equations (see, for example, \cite{Oksendal}) gives the above Equation (\ref{expected_exit}) for the expected exit time. For further details on the topic of this section, see section \ref{stochastics}.

\section{Analysis of response to a train-of-pulses}\label{section:train_of_pulses}
In addition to single-pulse experiments, 
TAP experimental studies often consider finely controlled multiple pulse conditions. Mathematically, it turns out that 
an infinite, periodic train of small pulses behaves well and  provides robust characteristics that can be compared with experimental data. 
Since gas exit flow is now periodic---with period $T$ given by the time interval between consecutive injected pulses---exit flow characteristics are naturally captured by Fourier coefficients, which can be effectively obtained from the single pulse exit flow function already derived, by means of the Poisson summation formula. This is shown in this final section.

Let us summarize here the main conclusions obtained in this section. Let $h(t)$ be the outlet flux generated by an
instantaneous pulse of amount $a$ at $x_0<L$ for the reaction-advection-diffusion problem with zero-flux boundary condition at $0$ and
absorbing  boundary condition at $L$. Extend $h$ to vanish for $t\leq 0$. For every $T>0$, the series
\[j_T(t):=\sum_{m\in \mathbb{Z}} h(t-mT)\]
defines a smooth (i.e., $C^\infty$) $T$-periodic function. Its Fourier coefficients are
\[\widehat{j}_{T,\ell}=\frac1TH(i\ell \Omega), \ \ \Omega=\frac{2\pi}{T},\]
where $H(s)$ is the Fourier-Laplace transform of $h(t)$:
\[H(s):=\int_0^\infty h(t) e^{-st}\, dt.\]
Moreover, 
\[\left| \widehat{j}_{T, \ell}\right|\leq C\exp\left(-c\sqrt{|\ell|}\right),\]
so its Fourier series and all differentiated Fourier series converge absolutely and uniformly. It will be further shown 
that $H(s)$ admits an   explicit description, which, for $x_0=0$ is:
\[H(s)=\frac{a e^{vL/(2D)}}{
\cosh\left(L\sqrt{\frac{v^2}{4D^2}+\frac{s+k}{D}}\right) +
 \frac{v}{2D}\frac{\sinh\left(L\sqrt{\frac{v^2}{4D^2}+\frac{s+k}{D}}\right)}{\sqrt{\frac{v^2}{4D^2}+\frac{s+k}{D}}}
}.\]
Finally, we show how the single-pulse response moments are encoded in the low-frequency behavior of the Fourier coefficients
of the periodic response case. Further observations are made along the way which we believe may be of
interest to the analysis of TAP experimental data.

\subsection{The equation for a periodic train of pulses} 
It may be helpful to summarize the notation used for the single-pulse problem. Equations will mostly be expressed
in the dimensional variables. The single-pulse problem is
\begin{equation}\label{main_equation_again}\frac{\partial c}{\partial t}=D\frac{\partial^2c}{\partial x^2}-v\frac{\partial c}{\partial x} -kc, \ \ 0<x<L,
\end{equation}
with boundary conditions $j(0,t)=0$, $c(L,t)=0$, where
$j(x,t)= vc(x,t)-D\frac{\partial c}{\partial x}(x,t),$
and an instantaneous pulse of amount $a$ at $x=x_0$ (typically taken to be very close to $0$). The measured response (exit flow) is
\[h(t):= j(L,t), \ \ \ t>0.\]
We have obtained the eigenfunction representation of the exit flow
\[h(t)= C\sum_{n=1}^\infty b_n e^{-\lambda_nt},\]
where 
\[C = \frac{aD}{L^2}\exp\left[ \frac{\text{Pe}}{2}\left(1-\frac{x_0}{L}\right)\right], \ b_n=\frac{\mu_n \sin\left(\mu_n \left(1-\frac{x_0}{L}\right)\right)}{\frac12 +\frac{\text{Pe}}{4}\left(\frac{\sin\mu_n}{\mu_n}\right)^2}, \ \lambda_n = k +\frac{\text{Pe}^2}{4t_d} +\frac{\mu_n^2}{t_d},\ t_d=\frac{L^2}{D}.\]

Let $T>0$ be the time between successive pulses and suppose that identical pulses of amount $a$ are injected at times
$t=mT$, $m=0, 1, 2, \dots$  and position $x_0$. One formulation of the problem is to solve a succession of initial-boundary value problems
so that, between two consecutive pulses Equation (\ref{main_equation_again}) holds, and, at the instant $t=mT$, the
concentration is updated according to 
\begin{equation}\label{discontinuous_condition}c(x,mT^+)=c(x,mT^-) + a \delta_{x_0}(x),\end{equation}
where $\delta_{x_0}(x)$ is the Dirac point mass distribution at $x_0$.  It is easily shown (by restating the problem in integral form) that
we may, equivalently, incorporate the pulses directly into the evolution equation  (\ref{main_equation_again}), and write
\begin{equation}
  \frac{\partial c}{\partial t}=D\frac{\partial^2c}{\partial x^2}-v\frac{\partial c}{\partial x} -kc+ a\delta_{x_0}(x)\sum_{m=0}^\infty \delta(t-mT).
\end{equation}
Across each pulse time, Equation (\ref{discontinuous_condition}) holds. Since the problem is linear, the exit flow generated by the
finite train is 
\[j_N(t) = \sum_{0\leq m\leq N, mT<t} h(t-mT).\]
Now write, after the $N$-th pulse, $t=NT +\theta$, $0<\theta <T$, and 
$j_N(NT+\theta)=\sum_{r=0}^N h(\theta+rT).$
The limiting periodic response is therefore
\[j_{\text{\tiny per}}(\theta) = \sum_{r=0}^\infty h(\theta +r T), \ \ 0\leq \theta<T.\]
Equivalently, after extending $h(t)$ to be zero for $t<0$, 
\begin{equation}\label{sum_h_equation}
j_{\text{\tiny per}}(\theta) = \sum_{m\in \mathbb{Z}}  h(\theta-mT).
\end{equation}
Our goal is to extract useful information about the process from $j_{\text{\tiny per}}(\theta) $.
\subsection{Convergence to the periodic exit flow}
\begin{figure}[htbp]
\begin{center}
\includegraphics[width=4.0in]{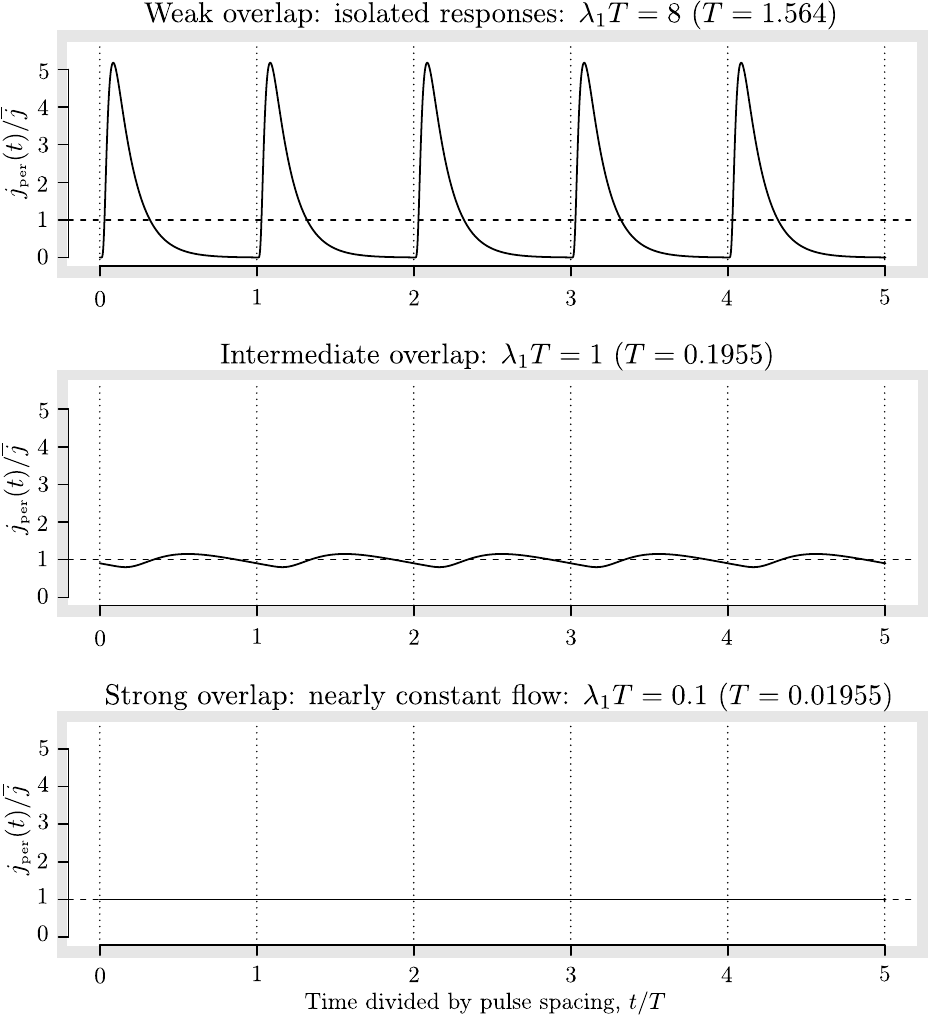}\ \ 
\caption{{\small Three regimes of gas exit flow defined in terms of the parameter $\lambda_1 T$.}}
\label{periodic_pulses}
\end{center}
\end{figure}

Since the operator has a spectral gap, the single pulse response decays exponentially. The slowest decay is 
\[\lambda_1 = k+\frac{\text{Pe}^2}{4t_d} +\frac{\mu_1^2}{t_d}>0. \]
Consequently, for sufficiently large $t$ we have $|h(t)| \leq C_1 e^{-\lambda_1 t}. $ It follows that
\[\left|j_{\text{\tiny per}}(\theta)-j_N(NT+\theta)\right|\leq C_2 \frac{e^{-\lambda_1\left(\theta + (N+1)T\right)}}{1-e^{-\lambda_1T}}.\]
Thus convergence from cycle to cycle is geometric, with approximate contraction factor $e^{-\lambda_1 T}$. This produces three useful
regimes:
\begin{enumerate}
\item $\lambda_1 T\gg 1$: successive pulses in the outflow are nearly isolated;
\item $\lambda_1 T\approx 1$: the responses to nearby pulses overlap appreciably;
\item $\lambda_1 T\ll 1$: strong accumulation and an almost constant outlet flow.  
\end{enumerate}
These regimes are illustrated in Figure \ref{periodic_pulses}.

\vspace{0.1in}
The quantity $e^{-\lambda_1T}$ is therefore a natural measure of pulse-to-pulse memory.

\subsection{Series description of the periodic response}\label{convergence_issues}
For $0<\theta<T$, the pulse number and eigenfunction sums  can be interchanged because the exponential factor $e^{-\lambda_n\theta}$ regularizes the eigenfunction expansion. Hence
\begin{align*}j_{\text{\tiny per}}(\theta) &= C\sum_{r=0}^\infty\sum_{n=1}^\infty b_n e^{-\lambda_n (\theta + rT)}\\
&=  C\sum_{n=1}^\infty b_n e^{-\lambda_n \theta}  \sum_{r=0}^\infty e^{-\lambda_n rT}\\
&= C\sum_{n=1}^\infty\frac{b_n e^{-\lambda_n \theta}}{1-e^{-\lambda_nT}},
\end{align*}
valid for $0<\theta<T$. It is important to note that this formula should not be evaluated at $\theta=0$. The correct endpoint value is instead
\[j_{\text{\tiny per}}(0) =\sum_{r=1}^\infty h(rT).\]
The contribution $r=0$ is $h(0)=0$, but its value cannot be recovered by naively setting $t=0$ in the eigenfunction series.

It is worth stopping at this point to consider convergence issues.  The eigenfunction series  for $h(t)$ is not absolutely or uniformly convergent down to $t=0$.  However, the exit flow  function itself is much better behaved. Let $d=L-x_0>0$ be the distance from the injection point to the outlet. Standard  heat-kernel behavior gives, for small $t$,
\[h(t)=O\left(t^{-3/2} \exp\left[-\frac{d^2}{4Dt}\right]\right).\]
More generally, for every $r\geq 0$, the derivatives of $h(t)$ satisfy
\[h^{(r)}(t)=O\left(t^{-N_r} \exp\left[-\frac{d^2}{4Dt}\right]\right)\]
for some $N_r$. Hence the limit of $h^{(r)}(t)$ as $t$ decreases to $0$ is $0$ for every $r$. The causal extension $h(t)=0$, $t\leq 0$, is therefore $C^\infty$ at $t=0$.  So the singular behavior is a defect of the termwise eigenfunction representation at the initial time, not a singularity of the physical  exit flow. Very high spatial modes undergo strong cancellation in order to produce the exponentially small short-time signal. This observation is central to applying the Poisson summation, which we do in the next subsection.

\subsection{Poisson summation and the Fourier coefficients}
Let $\Omega:=\frac{2\pi}{T}$. Write the Fourier series of the periodic exit flow as
\[j_{\text{\tiny per}}(t)=\sum_{\ell\in \mathbb{Z}}\widehat{j}_\ell e^{i\ell \Omega t},\]
where 
\[\widehat{j}_\ell =\frac1T \int_0^T j_{\text{\tiny per}}(t) e^{-i\ell\Omega t}\, dt.\]
Rewriting Equation(\ref{sum_h_equation}) as
$j_{\text{\tiny per}}(t)=\sum_{r=0}^\infty h(t+rT)$, $0\leq t<T$, 
we obtain
\[\widehat{j}_\ell =\frac1T\sum_{r=0}^\infty \int_0^T h(t+rT) e^{-i\ell \Omega t}\, dt=\frac1T \sum_{r=0}^\infty \int_{rT}^{(r+1)T} h(u) e^{-i\ell \Omega(u-rT)}\, du.\]
But $e^{i\ell \Omega r T}=e^{2\pi i\ell r}=1$. Therefore
\begin{equation}
\widehat{j}_\ell = \frac1T\int_0^\infty h(u) e^{-i\ell \Omega u}\, du=\frac1T H(i\ell\Omega),
\end{equation}
where we have introduced the one-sided Fourier-Laplace transform
\[H(s):= \int_0^\infty h(t) e^{-st}\, dt.\]
Consequently, we obtain the Poisson summation formula for the causal single-pulse response:
\begin{equation*}
j_{\text{\tiny per}}(t)=\frac1T\sum_{\ell\in \mathbb{Z}} H(i\ell \Omega) e^{i\ell \Omega t}.
\end{equation*}
\subsection{First arrival time and closed form of $H$}\label{closed_form_H}
For more details on the material of this section, see the final section \ref{stochastics}.

It is possible to obtain the transform $H(s)$ without summing the spatial eigenfunctions. Let us interpret $h(t)/a$ as the density of arrival at $x=L$, with
exponential killing at rate $k$. Its Laplace transform is 
\begin{equation}\label{Laplace_H}\frac1a H(s) = \mathbb{E}_{x_0}\left[e^{-(s+k)\tau_L}\right],\end{equation}
where $\tau_L$ is the first arrival time at $L$ for a diffusion with drift $v$, reflected at $0$. (See the explanations in Section \ref{stochastics}.)
The function \[u_s(x) = \mathbb{E}_x\left[e^{-(s+k)\tau_L}\right]\]
satisfies the backward boundary-value problem for the Feynman-Kac equation
\begin{equation}\label{feynman_kac}Du''_s  +vu_s' -(s+k)u_s =0, \ \ u_s'(0)=0, \  \ u_s(L)=1.\end{equation}
Introduce
\[\beta:=\frac{v}{2D}, \ \ \ q(s)=\sqrt{\beta^2 +\frac{s+k}{D}},\]
with  the square root chosen so that the real part of $q(s)$ is positive when the real part of $s$ is non-negative. 
Solving the ordinary differential equation (\ref{feynman_kac}), yields
\[u_s(x) = e^{\beta(L-x)}\frac{\cosh(qx) +\frac{\beta}{q}\sinh(qx)}{\cosh(qL) +\frac{\beta}{q}\sinh(qL)}.\]
We thus obtain the explicit form of $H(s)$:
\begin{equation*}
H(s)=ae^{\beta(L-x_0)}\frac{\cosh(q(s)x_0) +\frac{\beta}{q(s)}\sinh(q(s)x_0)}{\cosh(q(s)L) +\frac{\beta}{q(s)}\sinh(q(s)L)}.
\end{equation*}
For injection at $x_0=0$, which is the case in standard TAP experiments, this expression simplifies substantially:
\[H(s)=\frac{a e^{vL/(2D)}}{
\cosh\left(L\sqrt{\frac{v^2}{4D^2}+\frac{s+k}{D}}\right) +
 \frac{v}{2D}\frac{\sinh\left(L\sqrt{\frac{v^2}{4D^2}+\frac{s+k}{D}}\right)}{\sqrt{\frac{v^2}{4D^2}+\frac{s+k}{D}}}
}.\]
The Fourier coefficients are now obtained from this expression by setting $s=i\ell \Omega$.
\subsection{Non-dimensional form of the Fourier coefficients}
Let us recall the  coefficients 
\[\xi_0=\frac{x_0}{L}, \ \ \text{Pe}=\frac{vL}{D}, \ \ \kappa_d =k t_d, \ \ t_d = \frac{L^2}{D},\]
where $\xi_0$, $\text{Pe}$ and $\kappa_d$ are non-dimensional.  For a complex non-dimensional frequency $z=st_d$, we define
\[Q(z) =\sqrt{\frac{\text{Pe}^2}{4} +\kappa_d +z}.\]
Then, 
\[H(s)=H(z/t_d)= ae^{\frac{\text{Pe}}{2}(1-\xi_0)}
\frac{\cosh(Q(z)\xi_0) +\frac{\text{Pe}}{2Q(z)}\sinh(Q(z)\xi_0)}{
\cosh Q(z)+\frac{\text{Pe}}{2Q(z)}\sinh Q(z). 
} \]
For the $\ell$-th harmonic, $z=i\ell\Omega t_d =i\frac{2\pi \ell t_d}{T}.$ Thus
\[\widehat{j}_\ell =\frac{a}Te^{\frac{\text{Pe}}{2}(1-\xi_0)} 
\frac{\cosh(Q_\ell\xi_0) +\frac{\text{Pe}}{2Q_\ell}\sinh(Q_\ell\xi_0)}{
\cosh Q_\ell+\frac{\text{Pe}}{2Q_\ell}\sinh Q_\ell}, 
\]
where \[Q_\ell=\sqrt{\frac{\text{Pe}^2}{4} + \kappa_d + i\frac{2\pi \ell t_d}{T}}.\]
Thus the entire response is controlled by  the dimensionless parameters $\text{Pe}$, $\kappa_d$, $T/t_d$, $\xi_0$. 
\subsection{Absolute convergence of the Fourier series}
The closed form obtained for the Fourier coefficients resolves the convergence issue remarked earlier in section \ref{convergence_issues}.
For large $|\omega|$, 
\[q(i\omega) = (1+i\, \text{sgn}\,\omega) \sqrt\frac{|\omega|}{2D} +O(|\omega|^{-1/2}).\]
Consequently, 
\[|H(i\omega)|\leq C\exp\left[-(L-x_0)\sqrt{\frac{|\omega|}{2D}}\right].\]
At the harmonic frequencies $\omega=\ell \Omega$,
\[
|\widehat{j}_\ell| \leq \frac{C}{T} 
 \exp\left[
 -(L-x_0)\sqrt{\frac{\pi|\ell|}{DT}}
 \right], 
\]
so that, for all $r\geq 0$, 
$\sum_{\ell\in \mathbb{Z}}|\ell|^r |\widehat{j}_\ell|<\infty.$ This implies that $j_{\text{\tiny per}}$ is a smooth (i.e., $C^\infty$) function
on the circle $\mathbb{R}/T\mathbb{Z}$. It is thus interesting to contrast the convergence of the spatial eigenfunction series representation at the time of the injected pulse of gas, which is not absolutely convergent, with  the temporal Fourier series of the periodic outlet response, which is absolutely  and uniformly convergent. The spatial separation between injection and detection results in the strong high-frequency damping responsible for the better convergence in the latter case. 

\subsection{The zero-frequency coefficient}
The zero-frequency coefficient is of particular interest as it gives the average exit rate, which we denote $\overline{j}$:
\[\overline{j}=\widehat{j}_0=\frac{H(0)}{T}= \frac{a}T \mathbb{E}_{x_0}\left[e^{-k\tau_L}\right].\]
For $k=0$, all the gas injected in one pulse eventually exits, and $H(0)=a$. So $\overline{j}=a/T$ if $k=0$.
The mean outlet flow in the inert case is independent of $D$ and $v$, naturally. Diffusion and advection affect the oscillatory terms but not
the zero order (direct current) throughput.  With reaction, $H(0)/a<1$ and the deficit
$1-\frac{H(0)}{a}$ is the fraction consumed before reaching the outlet.

\subsection{Harmonic amplitudes and phases as experimental data}
Since $j_{\text{\tiny per}}$ is real, we have $\widehat{j}_{-\ell}=\overline{\widehat{j}_{\ell}}$. Thus
\[j_{\text{\tiny per}}=\widehat{j}_0+2\sum_{\ell=1}^\infty |\widehat{j}_\ell| \cos\left(\ell\Omega t +\phi_\ell\right),\]
where $\phi_\ell$ is the argument of $\widehat{j}_\ell$. 
The experimentally useful quantities are therefore
$\frac{|\widehat{j}_\ell|}{\widehat{j}_0}$ and $\phi_\ell$. The normalization by $\widehat{j}_0$ makes these quantities independent of the pulse amount $a$. They characterize the shape of the periodic response rather than its overall scale. In an experimental work,
a parameter-estimation scheme could be used to fit 
\[\left\{
\frac{|\widehat{j}_\ell|}{\widehat{j}_0}, \phi_\ell\right\}_{\ell=1}^N\]
to the theoretical expressions as functions of $D, v, k$.

\subsection{Connection with the single-pulse response moments}
Recall the raw moments
\[M_r = \int_0^\infty t^r h(t)\, dt\]
for the single-pulse response. They are directly encoded in the low-frequency behavior of the Fourier coefficients. In fact, since
\[H(s)=\int_0^\infty h(t) e^{-st}\,dt,\]
we have
$M_r=(-1)^rH^{(r)}(0).$
Hence
\[H(i\omega) = M_0 -i\omega M_1 -\frac{\omega^2}{2}M_2  +\frac{i\omega^3}{6}M_3+\cdots.\]
Let \[p(t):= \frac{h(t)}{M_0}\]
be the normalized exit-time distribution, and write
\[\overline{t}=\frac{M_1}{M_0}, \ \ \ \sigma_t^2=\frac{M_2}{M_0}- \left(\frac{M_1}{M_0}\right)^2.\]
Then
\[\frac{H(i\omega)}{H(0)}=1-i\omega\overline{t} -\frac{\omega^2}{2}\frac{M_2}{M_0}+O(\omega^3),\]
and
\[\log\frac{H(i\omega)}{H(0)} = -i\omega \overline{t} -\frac{\omega^2}{2}\sigma_t^2 +\frac{i\omega^3}{6}\kappa_3 + O(\omega^4),\]
where $\kappa_3$ is the third cumulant of $p(t)$. Therefore, for the first few harmonics when $\Omega$ is small,
\[\phi_\ell \approx -\ell \Omega \overline{t}, \ \ \ \log\frac{|\widehat{j}_\ell|}{\widehat{j}_0}\approx -\frac12 (\ell \Omega)^2\sigma_t^2.\]
Thus, we have  the following observations when $\Omega$ is small (equivalently,  $T$ is large),
which provide a direct connection between the analysis of  moments for a single pulse, developed  earlier in the paper, to 
the periodic-pulse problem:
\begin{enumerate}
\item The phase measures the mean exit time.
\item The quadratic low-frequency decay of the log-amplitudes measures the exit-time variance.
\item Departures from these approximations encode skewness and higher cumulants.
\end{enumerate}

\section{Stochastic calculus perspective}\label{stochastics}
We provide here an overview of our analysis from an It\^o calculus perspective. 
\subsection{The SDE for molecular transport}
Let us first consider the process without reaction, $k=0$, and consider the stochastic  motion of a single molecule. The diffusion-advection transport equation is the Fokker-Planck equation for the molecular motion $X_t$ in the interior interval $(0,L)$, which satisfies
the stochastic differential equation
\[dX_t = v\, dt +\sqrt{2 D} \, dW_t,\]
where $W_t$ is a standard Wiener process.  The infinitesimal {\em backward generator} of this diffusion is 
\[\mathcal{L} f = Df'' + vf'.\]
On the other hand, its probability density $p(y,t|x)$, conditional on starting at $x$, satisfies the {\em forward}, or {Fokker-Planck} equation
\[\frac{\partial p}{\partial t}=\mathcal{L}^* p = D\frac{\partial^2p}{\partial y^2} -v\frac{\partial p}{\partial y}.\]
This is the transport part of the concentration equation (\ref{RAD_equation}), with flux $j = vc -D\frac{\partial c}{\partial x}$, zero flux at $x=0$ and
$c(L,t)=0$. Note that $\mathcal{L}$ acts on functions of the starting position, whereas $\mathcal{L}^*$ acts on probability densities. 
The boundary conditions correspond to Brownian motion with drift {\em reflected} at $0$ and stopped at $L$.
One may write it as the reflected SDE
\[dX_t = v\, dt + \sqrt{2D} \, dW_t +dK_t,\]
where $K_t$ is nondecreasing finite-variation process that increases only   when $X_t=0$;
in one dimension it is related, according to normalization convention, to the boundary local time. For the forward density, reflection means zero probability current, which gives the Robin condition
\[vp(0,t|x)-D\frac{\partial p}{\partial y}(0,t|x)=0,\]
but for the backward generator, reflection appears as the Neumann condition $f'(0)=0$.  
  
 A key observation is that the exit flow is the first-passage time density. More precisely, 
 let \[\tau_L :=\inf \{t\geq 0: X_t=L\}\]
be the time at which the molecule first exits the reactor. If it begins at $x_0$, then
\[S(t)=\text{Pr}(\tau_L>t)=\int_0^L p(x,t|x_0)\, dx\]
is its survival probability inside the reactor. Integrating the Fokker-Planck equation gives
\[\frac{dS}{dt} = j(0,t)-j(L,t) = -j(L,t),\]
for a unit-mass pulse. Note the use of the boundary condition $j(0,t)=0$. But \[f_{\tau_L}(t)=-\frac{d}{dt}\text{Pr}(\tau_L>t)\]
is the probability density of the first-passage time.  Therefore
\[\frac1a j_0(L,t)=f_{\tau_L}(t).\]
In words, the normalized inert-gas exit flow curve is precisely the probability density of the first-passage time at $L$ of a reflected diffusion through the reactor. 
\subsection{Reaction as independent exponential killing}
Let us now restore the reaction term $-kc$. A first-order irreversible reaction has a clear stochastic interpretation: give each molecule
an independent exponential lifetime $\zeta\sim \text{Exp}(k)$.
A molecule is detected at the outlet only if $\tau_L<\zeta.$ Conditional on $\tau_L=t$,
\[\text{Pr}(\zeta>t)=e^{-kt}.\]
Consequently,
\[j_k(L,t)=ae^{-kt}f_{\tau_L}(t)=e^{-kt}j_0(L,t).\]
This is the probabilistic interpretation of the previously analytically  derived (from the series solution) Equation (\ref{ratio}).
To emphasize the point, we can write the normalized exit flow in the presence of the irreversible first-order reaction with constant $k$ as
the product of a first-passage-time density and a survival probability:
\[\frac1a j_k(L,t)=\underbrace{f_{\tau_L}(t)}_{\text{molecule arrives at outlet at time t}} \times  \underbrace{e^{-kt}}_\text{molecule survives reaction until time $t$}.\]

\subsection{The  function $H(s)$}
Recall that  $H(s)=\int_0^\infty e^{-st}j_k(L,t)\, dt$. Using the previously obtained factorization identity, we arrive at
\[\frac{H(s)}{a} =\int_0^\infty e^{-st}e^{-kt} f_{\tau_L}(t)\, dt=\int_0^\infty e^{-(s+k)t} f_{\tau_L}(t)\, dt = \mathbb{E}_{x_0}\left[e^{-(s+k)\tau_L}\right].\]
This is the justification for Equation (\ref{Laplace_H}). In particular,
\[\frac{H(0)}a=\mathbb{E}_{x_0}\left[e^{-k\tau_L}\right]=\text{Pr}\left(\tau_L<\zeta\right). \]
This is the fraction of (unreacted) gas that eventually exits.
\subsection{The function $u_s(x)$ satisfies the backward equation}
 Set $q=s+k$ and  define
 \[u_q(x)=\mathbb{E}_x\left[e^{-q\tau_L}\right].\]
 The determination of $u_q(x)$ is a Feynman-Kac problem, although in this particularly simple situation one can derive the equation directly
 from Dynkin's formula or It\^o's formula. The result is
 \[Du_q''(x) +vu_q'(x) -qu_q(x)=0, \ \ 0<x<L,\]
 with boundary conditions
 \[u'_q(0)=0, \ \ \ u_q(L)=1. \]
 Perhaps a brief explanation of this boundary-value problem is justified. The equation can be written as $(\mathcal{L}-q)u_q=0$, in which $q$ is the killing rate that accounts for the exponential term $e^{-qt}$. The reflecting boundary gives $u'_q(0)=0$, and if the molecule starts at the outlet, then $\tau_L=0$ so $u_q(L) = 1$. The second order equation for $u_q(x)$ can be obtained as follows. Consider the process
 \[M_t = e^{-q(t\wedge \tau_L)} u_q\left(X_{t\wedge \tau_L}\right).\]
Applying It\^o's formula to the reflected process produces,  
\[dM_t = e^{-qt}\left[(\mathcal{L}u_q-qu_q)\, dt + \sqrt{2D} u_q' (X_t)\, dW_t + u_q'(0) \, dK_t\right].\]
If $(\mathcal{L}-q)u_q=0$ and $u_q'(0)=0$, the drift and reflection terms disappear, leaving a martingale. 
Stopping at $\tau_L$ and taking expectations gives
\[u_q(x)=\mathbb{E}_x\left[e^{-q\tau_L} u_q(X_{\tau_L})\right]=\mathbb{E}_x\left[e^{-q\tau_L} \right],\]
where we have used that $X_{\tau_L}=L$ and $u_q(L)=1$. This is the justification of the backward boundary-value problem used in section \ref{closed_form_H}. The explicit solution given in that section is obtained by elementary methods for solving linear ordinary differential equations with constant coefficients.
\subsection{Further remarks concerning section \ref{sec05}}
In section \ref{sec05} we introduced $T(x)=\mathbb{E}_x[\tau_L]$ and stated that Dynkin's formula gives $DT''+vT'=-1$ with boundary conditions
$T'(0)=0$, $T(L)=0$. This equation was then solved to obtain the mean exit time and its diffusion- and advection-dominated limits.  
Note that 
\[T(x)=-\left.\frac{\partial u_q(x)}{\partial q}\right|_{q=0}.\]
Since $u_0(x)=1$, we may differentiate $\mathcal{L}u_q -q u_q=0$ with respect to $q$ at zero to obtain
\[\mathcal{L}\left.\frac{\partial u_q}{\partial q}\right|_{q=0}= 1.\]
But this is precisely the above equation for $T(x)$. The boundary values are similarly obtained by differentiating in $q$ the boundary conditions for $u_q(x)$. In fact, all first-passage moments fit into this same scheme.  Define
\[m_r(x):=\mathbb{E}_x\left[\tau_L^r\right].\]
Since 
\[u_q(x)=\mathbb{E}_x\left[e^{-q\tau_L}\right]=\sum_{r=0}^\infty\frac{(-q)^r}{r!}m_r(x)\]
near $q=0$, differentiation of the backward equation gives 
$\mathcal{L} m_r =-r m_{r-1},$
or \[Dm_r'' +vm_r'=-r m_{r-1},\]
a recursion equation for the moments,  with the boundary conditions  $m_r'(0)=0$,  $m_r(L)=0$ for $r\geq 1$ and $m_0=1$. 
This provides a general probabilistic interpretation of the moments. For an inert gas,
$M_r=a \mathbb{E}_{x_0}\left[\tau_L^r\right].$ With reaction $M_r=a \mathbb{E}_{x_0}\left[\tau_L^r e^{-k\tau_L}\right].$
In particular, 
\[\frac{M_1}{M_0} = \frac{ \mathbb{E}_{x_0}\left[\tau_L e^{-k\tau_L}\right]}{ \mathbb{E}_{x_0}\left[e^{-k\tau_L}\right]} = \mathbb{E}_{x_0}\left[\tau_L|\tau_L<\zeta\right].\]
This is the mean residence time among the molecules that survive reaction and actually reach the detector.
Note that, for $k=0$, we have $M_1/M_0=\mathbb{E}_{x_0}[\tau_L].$ This is  the relation $t_\text{\tiny moments}=t_{\text{\tiny mean}}$  claimed in section \ref{sec05}.
\subsection{Connections to the periodic-pulse problem}
For the periodic pulse train it was found that 
$\widehat{j}_\ell=\frac1T H(i\ell \Omega)$, where  $\Omega=\frac{2\pi}{T}.$
The stochastic representation  gives
\[\widehat{j}_\ell = \frac{a}T \mathbb{E}_{x_0}\left[e^{-k\tau_L }e^{-i\ell\Omega\tau_L}\right].\]
This provides the Fourier coefficients with a very concrete interpretation: {\em they are samples of the Fourier transform of
the first passage time  distribution, exponentially weighted by the probability of chemical survival.} For $k=0$,
\[\frac{T}{a} \widehat{j}_\ell = \mathbb{E}_{x_0}\left[e^{-i\ell \Omega \tau_L}\right]\]
is  simply the {\em characteristic function of the residence time}, sampled at the harmonic frequencies.

\section{Conclusions}
A mathematical framework for the description  of pulse-response studies of gas-solid catalytic reactions  was presented based on an advection-diffusion-reaction model. 
For the transport model consisting of Fickian diffusion plus advection, the transport curve was found and some of its properties, such as
moments of the exit flow and peak characteristics, were described analytically and numerically.  The same was obtained for the process with a first-order irreversible reaction 
added. We analyze both the transient response of a single pulse system and the periodic response to a train of pulses via Poisson summation,
obtaining explicit expressions  for the Fourier coefficients of the periodic exit flow.
We thus obtain  numerical characteristics of the system that should prove  useful in the analysis of experimental data.
The mathematical analysis is approached from a  partial differential equations and a stochastic calculus perspective. 
The obtained results strengthen  the theoretical foundations for pulse-response studies, opening new perspectives for extracting transport and kinetic dependencies of heterogeneous catalytic processes in the wide domain.

\end{document}